\documentclass[12pt]{article}
\usepackage{amssymb}
\usepackage{amsthm}
\usepackage{amsmath}
\usepackage{graphicx}
 \newtheorem{thm}{Theorem}[section]

 \theoremstyle{definition}
 
 \newtheorem{rem}[thm]{Remark}
 \numberwithin{equation}{section}

\begin{document}
\title{Logarithmically Improved Serrin's Criteria for Navier-Stokes Equations}
\author{Yi Zhou\footnote{School of Mathematical Sciences, Fudan
  University, Shanghai 200433, China. {\it Email: yizhou@fudan.ac.cn}} \and
  Zhen Lei\footnote{School of Mathematical Sciences, Fudan
  University, Shanghai 200433, China. {\it Email:
  leizhn@yahoo.com}}}
\date{\today}
\maketitle
\begin{abstract}
In this paper we prove the logarithmically improved Serrin's
criteria to the three-dimensional incompressible Navier-Stokes
equations.
\end{abstract}

\section{Introduction}%----------------------------introduction

As is well-known that the three-dimensional incompressible
Navier-Stokes equations in $\mathbb{R}^3$ take the form
\begin{equation}\label{NS}
\begin{cases}
u_t + u\cdot\nabla u + \nabla q
  = \mu\Delta u,\\[-4mm]\\
\nabla\cdot u = 0,
\end{cases}
\end{equation}
where $u = (u_1, u_2, u_3)^T$ is the velocity of the flows, $q$ is
the scalar pressure and $\mu$ is the viscosity of the fluid which
is the inverse of the Reynolds number. We refer the reader to
\cite{Lady70, Temam, CF, MB} for general descriptions of Euler and
Navier-Stokes equations.

Despite a great deal of efforts by mathematicians and physicists,
the question of whether a solution of the 3D incompressible
Navier-Stokes equations can develop a finite time singularity from
smooth initial data with finite energy is still one of the most
outstanding mathematical open problems \cite{Fefferman}. In the
absence of a well-posedness theory, the development of
blowup/non-blowup criteria is of major importance for both
theoretical and practical purposes. There has been a lot of
progress on Euler and Navier-Stokes equations along this
direction. For example, by the well-known Serrin's criteria
\cite{Prodi, Serrin1, Serrin2, Struwe} for Navier-Stokes
equations, any Leray-Hopf weak solution $u$ is smooth for $0 \leq
t \leq T$ provided that
\begin{equation}\nonumber
\int_0^T\|u(t, \cdot)\|_{L^p}^s\ dt < \infty
\end{equation}
holds for any pair of constants $(p, s)$ with $\frac{3}{p} +
\frac{2}{s} \leq 1$ and $3 < p \leq \infty$. In fact, the highly
nontrivial case of $p = 3$ is also true, which was only proved
recently by Iskauriaza, Seregin and Sverak
\cite{IskauriazaSereginShverak}. By the celebrated
Beale-Kato-Majda's criterion \cite{BKM84}, if
\begin{equation}\nonumber
\int_0^T\|\omega(t, \cdot)\|_{L^\infty}\ dt < \infty,
\end{equation}
where $\omega = \nabla\times u$ is the vorticity, then $u$ is
smooth at time $T$. The above Beale-Kato-Majda's criterion was
slightly improved by Kozono-Taniuchi \cite{KozonoTaniuchi} by
replacing $\|\omega(t, \cdot)\|_{L^\infty}$ as $\|\omega(t,
\cdot)\|_{\rm BMO}$. Let us also cite a result by Constantin and
Fefferman \cite{CFefferman} who gave a condition involving only
the direction of the vorticity and a condition involving the lower
bound of the pressure by Seregin and Sverak \cite{SereginSverak}.

Our purpose of this paper is to establish the logarithmically
improved Serrin's criterion for three-dimensional Navier-Stokes
equations.
\begin{thm}\label{thmLogSerrin}
Suppose that $u_0 \in H^2(\mathbb{R}^3)$ and $u$ is a smooth
solution to the incompressible Navier-Stokes equations \eqref{NS}
with the initial data $u_0(x)$ for $0 \leq t < T$. Then $u$ is
smooth at time $t = T$ provided that
\begin{equation}\label{LogSerrin}
\int_0^T\frac{\|u(t, \cdot)\|_{L^p}^s}{1 + \ln\big(e + \|u(t,
\cdot)\|_{L^\infty}\big)}\ dt < \infty
\end{equation}
for any pair of constants $(p, s)$ with $\frac{3}{p} + \frac{2}{s}
\leq 1$ and $3 < p \leq \infty$.
\end{thm}

\begin{rem}
Serrin's criteria were first obtained by Prodi \cite{Prodi} and
Serrin \cite{Serrin1} in the whole space case. The local version
was established by Serrin \cite{Serrin2} when $\frac{3}{p} +
\frac{2}{s} < 1$ and Struwe when $\frac{3}{p} + \frac{2}{s} = 1$.
For simplicity, we just focus on the whole space case and do not
pursue the bounded domain case. Moreover, in Theorem
\ref{LogSerrin}, we just assume that $u$ is a smooth solution for
$0 \leq t < T$. However, the conclusion in Theorem \ref{LogSerrin}
is still true for Leray-Hopf weak solutions, which can be proved
by the same smoothing technique as in Struwe \cite{Struwe}.
\end{rem}

\begin{rem}
We remark here that recently Chan and Vasseur \cite{ChanVasseur}
proved a logarithmically improved Serrin's criteria under assuming
that
\begin{equation}\nonumber
\int_0^T\int_{\mathbb{R}^3}\frac{|u(t, x)|^5}{\ln\big(e + |u(t,
x)|\big)}\ dxdt < \infty.
\end{equation}
Noting that
\begin{equation}\nonumber
\int_0^T\frac{\|u(t, \cdot)\|_{L^5}^5}{1 + \ln\big(e + \|u(t,
\cdot)\|_{L^\infty}\big)}\ dt \leq
\int_0^T\int_{\mathbb{R}^3}\frac{|u(t, x)|^5}{\ln\big(e + |u(t,
x)|\big)}\ dxdt,
\end{equation}
we find that our result covers the result in \cite{ChanVasseur} by
letting $p = s = 5$ in \eqref{LogSerrin}. Moreover, Chan and
Vasseur use De Giorgi's method, while our proof is just based on
energy method and is much simpler. After the completion of the
paper, we learned that the authors in \cite{ZhouY} have
independently proved results similar to those in Theorem
\ref{thmLogSerrin} in the framework of multiplier space.
\end{rem}

%\begin{rem}
%We also remark that recently Chen, Strain, Tsai and Yau
%\cite{CSTY07a, CSTY07b},  Koch, Nadirashvili, Seregin and Sverak
%\cite{KNSS}, Seregin andSverak \cite{SereginSverak2} studied the
%lower bound of the blowup rate of axi-symmetric Navier-Stokes
%equations. In particular, they obtained the regularity of
%solutions at $(T, 0)$ provided that the velocity field satisfies
%\begin{equation}\label{a1}
%|u(t, x)| \leq \frac{C_\star}{\sqrt{T - t}}.
%\end{equation}
%In fact, it is also of course an important question to make clear
%whether the solutions of Navier-Stokes equations blow up at $(T,
%x_\star)$ provided that
%\begin{equation}\label{a2}
%|u(t, x)| \leq \frac{C_\star}{\sqrt{|x - x_\star|^2 + T - t}}\
%{\rm or\ even}\  \frac{C_\star}{\sqrt{T - t}}.
%\end{equation}
%is satisfied in general case. See more discussions in
%\cite{CSTY07a}. Let $u$ satisfy \eqref{a2} and choose $p =
%\infty$, $s = 2$ in \eqref{LogSerrin}. It is easy to check that
%\begin{equation}\nonumber
%\frac{\|u(t, \cdot)\|_{L^\infty}^2}{1 + \ln\big(e + \|u(t,
%\cdot)\|_{L^\infty}\big)}\ {\rm is\ an\ increasing\ function\ of\
%\|u(t, \cdot)\|_{L^\infty}}
%\end{equation}
%when $\|u(t, \cdot)\|_{L^\infty}$ is large. Consequently,
%\begin{equation}\nonumber
%\int_0^t\frac{\|u(\tau, \cdot)\|_{L^\infty}^2}{1 + \ln\big(e +
%\|u(\tau, \cdot)\|_{L^\infty}\big)}d\tau \leq C_\star^2\ln\big(1 +
%\ln\frac{1}{T - t}\big),
%\end{equation}
%while Serrin's criterion only gives that
%\begin{equation}\nonumber
%\int_0^t\|u(\tau, \cdot)\|_{L^\infty}^2d\tau \leq C_\star^2
%\ln\frac{T}{T - t}.
%\end{equation}
%\end{rem}

\section{Logarithmically Improved Serrin's Criterion
for Navier-Stokes Equations}

In this section, we establish the logarithmically improved
Serrin's criterion for Navier-Stokes Equations and prove Theorem
\ref{thmLogSerrin}.

First of all, for any smooth solution $u$ to the three-dimensional
Euler and Navier-Stokes equations, one has the well-known energy
law:
\begin{eqnarray}\label{c1}
\frac{1}{2}\frac{d}{dt}\|u\|_{L^2}^2 + \mu\|\nabla u\|_{L^2}^2 =
0.
\end{eqnarray}
Next, let us apply $\nabla^3$ to \eqref{NS} and then take the
$L^2$ inner product of the resulting equations with $\nabla^3u$.
The standard energy estimate gives
\begin{eqnarray}\label{c2}
\frac{1}{2}\frac{d}{dt}\|\nabla^3u\|_{L^2}^2 + \mu\|\nabla^4
u\|_{L^2}^2 = - \int_{\mathbb{R}^3}\nabla^3u\nabla^3qdx  -
\int_{\mathbb{R}^3}\nabla^3u\nabla^3(u\cdot\nabla u)dx.
\end{eqnarray}
Noting the incompressible constraint $\nabla\cdot u = 0$ and using
integration by parts, one has
\begin{eqnarray}\nonumber
&&\frac{1}{2}\frac{d}{dt}\|\nabla^2u\|_{L^2}^2 + \mu\|\nabla^3
  u\|_{L^2}^2\\\nonumber
&&= -  2\int_{\mathbb{R}^3}\nabla^2u(\nabla u\cdot\nabla
  \nabla u)dx -  \int_{\mathbb{R}^3}\nabla^2u\nabla^2u
  \cdot\nabla udx\\\nonumber
&&\leq 5\|u\|_{L^p}\|\nabla^2u\|_{L^{\frac{2p}{p -
  2}}}\|\nabla^3u\|_{L^2},
\end{eqnarray}
for $3 < p \leq \infty$, where we used H${\rm \ddot{o}}$lder
inequality and integration by parts in the last inequality. Since
$3 < p \leq \infty$, one has $2 \leq \frac{2p}{p - 2} < 6$.
Consequently, by the following standard multiplicative inequality
\begin{eqnarray}\nonumber
\|\nabla^2u\|_{L^{\frac{2p}{p - 2}}} \leq C\|\nabla^2 u\|_{L^2}^{1
- \frac{3}{p}}\|\nabla^3u\|_{L^2}^{\frac{3}{p}},
\end{eqnarray}
we have
\begin{eqnarray}\label{d1}
&&\frac{1}{2}\frac{d}{dt}\|\nabla^2u\|_{L^2}^2 + \mu\|\nabla^3
  u\|_{L^2}^2\\\nonumber
&&\leq C\|u\|_{L^p}\|\nabla^2 u\|_{L^2}^{1 - \frac{3}{p}}
  \|\nabla^3u\|_{L^2}^{1 + \frac{3}{p}}\\\nonumber
&&\leq C(p, \mu)\|u\|_{L^p}^{\frac{2p}{p - 3}}\|\nabla^2
  u\|_{L^2}^2 + \frac{\mu}{2}\|\nabla^3u\|_{L^2}^2\\\nonumber
&&\leq C(p, \mu)\frac{\|u\|_{L^p}^{\frac{2p}{p - 3}}}{1 +
  \ln\big(e + \|u\|_{L^\infty}\big)}\big[1 + \ln\big(e +
  \|u\|_{L^\infty}\big)\big]\|\nabla^2u\|_{L^2}^2
  + \frac{\mu}{2}\|\nabla^3u\|_{L^2}^2.
\end{eqnarray}
Noting \eqref{c1}, we derive from \eqref{d1} that
\begin{eqnarray}\label{d2}
&&\frac{d}{dt}\|\nabla^2u\|_{L^2}^2 + \mu\|\nabla^3
  u\|_{L^2}^2\\\nonumber
&&\leq 2C(p, \mu)\frac{\|u\|_{L^p}^{\frac{2p}{p - 3}}}{1 +
  \ln\big(e + \|u\|_{L^\infty}\big)}\big[1 + \ln\big(e +
  \|\nabla^2u\|_{L^2}^2\big)\big]\|\nabla^2u\|_{L^2}^2.
\end{eqnarray}
Then using Gronwall's inequality, we have
\begin{eqnarray}\label{d3}
&1 + \ln\big(e + \|\nabla^2u(t, \cdot)\|_{L^2}^2\big) \leq \big[1
  + \ln\big(e + \|\nabla^2u_0\|_{L^2}^2\big)\big]\\\nonumber
&\quad\quad\quad\times\exp\Big\{2C(p, \mu)\int_0^t\frac{\|u(\tau,
  \cdot)\|_{L^p}^{\frac{2p}{p - 3}}}{1 +
  \ln\big(e + \|u\|_{L^\infty}\big)}\ d\tau\Big\},
\end{eqnarray}
which gives a finite bound for $\|\nabla^2u(T, \cdot)\|_{L^2}^2$
provided that
\begin{eqnarray}\nonumber
\int_0^T\frac{\|u(\tau,
  \cdot)\|_{L^p}^{\frac{2p}{p - 3}}}{1 +
  \ln\big(e + \|u\|_{L^\infty}\big)}\ d\tau
\end{eqnarray}
is finite.

%----------------------------------------------------------------------------acknowledgement
\section*{Acknowledgments} The work was partially supported by
the National Science Foundation of China under grant 10225102 and
a 973 project of the National Sciential Foundation of China. Zhen
Lei was partially supported by the NSFC 10801029 and PSFC
20070410160.

%-----------------------------------------------------------------------------bibliography

\end{document}